# SEQUENTIAL CHANGE-POINT DETECTION WHEN UNKNOWN PARAMETERS ARE PRESENT IN THE PRE-CHANGE DISTRIBUTION[1]

### By Yajun Mei

*California Institute of Technology and Fred Hutchinson Cancer Research Center*

In the sequential change-point detection literature, most research specifies a required frequency of false alarms at a given pre-change distribution $f_\theta$ and tries to minimize the detection delay for every possible post-change distribution $g_\lambda$. In this paper, motivated by a number of practical examples, we first consider the reverse question by specifying a required detection delay at a given post-change distribution and trying to minimize the frequency of false alarms for every possible pre-change distribution $f_\theta$. We present asymptotically optimal procedures for one-parameter exponential families. Next, we develop a general theory for change-point problems when both the pre-change distribution $f_\theta$ and the post-change distribution $g_\lambda$ involve unknown parameters. We also apply our approach to the special case of detecting shifts in the mean of independent normal observations.

**1. Introduction.** Suppose there is a process that produces a sequence of independent observations $X_1, X_2, \ldots$. Initially the process is "in control" and the true distribution of the $X$'s is $f_\theta$ for some $\theta \in \Theta$. At some unknown time $\nu$, the process goes "out of control" in the sense that the distribution of $X_\nu, X_{\nu+1}, \ldots$ is $g_\lambda$ for some $\lambda \in \Lambda$. It is desirable to raise an alarm as soon as the process is out of control so that we can take appropriate action. This is known as a change-point problem, or quickest change detection problem. By analogy with hypothesis testing terminology [12], we will refer to $\Theta$ ($\Lambda$) as a "simple" pre-change (post-change) hypothesis if it contains a single point and as a "composite" pre-change (post-change) hypothesis if it contains more than one point.

Received February 2004; revised February 2005.
[1]Supported in part by NIH Grant R01 AI055343.

The change-point problem originally arose from statistical quality control, and now it has many other important applications, including reliability, fault detection, finance, signal detection, surveillance and security systems. Extensive research has been done in this field during the last few decades. For recent reviews, we refer readers to [1, 9] and the references therein.

In the simplest case where both $\Theta$ and $\Lambda$ are simple, that is, the pre-change distribution $f_\theta$ and the post-change distribution $g_\lambda$ are completely specified, the problem is well understood and has been solved under a variety of criteria. Some popular schemes are Shewhart's control charts, moving average control charts, Page's CUSUM procedure and the Shiryayev–Roberts procedure; see [1, 17, 24, 25, 26]. The first asymptotic theory, using a minimax approach, was provided in [14].

In practice, the assumption of known pre-change distribution $f_\theta$ and post-change distribution $g_\lambda$ is too restrictive. Motivated by applications in statistical quality control, the standard formulation of a more flexible model assumes that $\Theta$ is simple and $\Lambda$ is composite, that is, $f_\theta$ is completely specified and the post-change distribution $g_\lambda$ involves an unknown parameter $\lambda$. See, for example, [9, 10, 11, 14, 20, 21, 29]. When the true $\theta$ of the pre-change distribution $f_\theta$ is unknown, it is typical to assume that a training sample is available so that one can use the method of "point estimation" to obtain a value $\theta_0$. However, it is well known that the performances of such procedures are very sensitive to the error in estimating $\theta$; see, for example, [30]. Thus we need to study change-point problems for composite pre-change hypotheses, which allow a range of "acceptable" values of $\theta$.

There are a few papers in the literature that use a parametric approach to deal with the case when the pre-change distribution involves unknown parameters (see, e.g., [6, 8, 22, 33, 34]), but all assume the availability of a training sample and/or the existence of an invariant structure. In this paper, we make no such assumptions. Our approach is motivated by the following examples.

EXAMPLE 1.1 (Water quality).  Suppose we are interested in monitoring a contaminant, say antimony, in drinking water. Because of its potential health effects, the U.S. Environmental Protection Agency (EPA) sets a maximum contaminant level goal (MCLG) and a maximum contaminant level (MCL). An MCLG is a nonenforceable but desirable health-related goal established at the level where there is no known or expected risk to health. An MCL is the enforceable limit set as close to the MCLG as possible. For antimony, both MCL and MCLG are 0.006 mg/L. Thus the water quality is "in control" as long as the level of the contaminant is less than MCLG, and we should take prompt action if the level exceeds MCL.

EXAMPLE 1.2 (Public health surveillance). Consider the surveillance of the incidence of rare health events. If the underlying disease rate is greater than some specified level, we want to detect it quickly so as to enable early intervention from a public health point of view and to avoid a much greater tragedy. Otherwise, the disease is "in control."

EXAMPLE 1.3 (Change in variability). In statistical process control, sometimes one is concerned about possible changes in the variance. When the value of the variance is greater than some pre-specified constant, the process should be stopped and declared "out of control." However, when the process is in control, there typically is no unique target value for the variance, which should be as small as the process permits.

EXAMPLE 1.4 (Signal disappearance). Suppose that one is monitoring or tracking a weak signal in a noisy environment. If the signal disappears, one wants to detect the disappearance as quickly as possible. Parameters $\theta$ associated with the signal, for example, its strength, are described by a composite hypothesis before it disappears, but by a simple hypothesis (strength equal to zero) afterward.

The essential feature of these examples is that the need to take action in response to a change in a parameter $\theta$ can be defined by a fixed threshold value. This inspires us to study change-point problems where $\Theta$ is composite and $\Lambda$ is simple. Unlike the standard formulation which specifies a required frequency of false alarms, our formulation specifies a required detection delay and seeks to minimize the frequency of false alarms for all possible pre-change distributions $f_\theta$. Section 2 uses this formulation to study the problem of detecting a change of the parameter value in a one-parameter exponential family. It is worthwhile pointing out that the generalized likelihood ratio method does not provide asymptotically optimal procedures under our formulation.

It is natural to combine the standard formulation with our formulation by considering change-point problems when both $\Theta$ and $\Lambda$ are composite, that is, both the pre-change distribution and the post-change distribution involve unknown parameters. Ideally we want to optimize all possible false alarm rates and all possible detection delays. Unfortunately this cannot be done, and there is no attractive definition of optimality in the literature for this problem. In Section 3, we propose a useful definition of "asymptotically optimal to first order" procedures, thereby generalizing Lorden's asymptotic theory, and develop such procedures with the idea of "optimizer."

This paper is organized as follows. In the remainder of this section we provide some notation and definitions based on the classical results for the change-point problem when both $\Theta$ and $\Lambda$ are simple. Section 2 establishes the asymptotic optimality of our proposed procedures for the problem of detecting a change of the parameter value in a one-parameter exponential family, and Section 3 develops an asymptotic theory for change-point problems when both the pre-change distribution and the post-change distribution involve unknown parameters. Both Sections 2 and 3 contain some numerical simulations. Section 4 illustrates the application of our general theory to the problem of detecting shifts in the mean of independent normal observations. Section 5 contains the proof of Theorem 2.1.

Denote by $\mathbf{P}_{\theta,\lambda}^{(\nu)}, \mathbf{E}_{\theta,\lambda}^{(\nu)}$ the probability measure and expectation, respectively, when $X_1, \ldots, X_{\nu-1}$ are distributed according to a pre-change distribution $f_\theta$ for some $\theta \in \Theta$ and $X_\nu, X_{\nu+1}, \ldots$ are distributed according to a post-change distribution $g_\lambda$ for some $\lambda \in \Lambda$. We shall also use $\mathbf{P}_\theta$ and $\mathbf{E}_\theta$ to denote the probability measure and expectation, respectively, under which $X_1, X_2, \ldots$ are independent and identically distributed with density $f_\theta$ (corresponding to $\nu = \infty$). In change-point problems, a procedure for detecting that a change has occurred is defined as a stopping time $N$ with respect to $\{X_n\}_{n\geq 1}$. The interpretation of $N$ is that, when $N = n$, we stop at $n$ and declare that a change has occurred somewhere in the first $n$ observations. The performance of $N$ is evaluated by two criteria: the long and short average run lengths (ARL). The long ARL is defined by $\mathbf{E}_\theta N$. Imagining repeated applications of such procedures, practitioners refer to the frequency of false alarms as $1/\mathbf{E}_\theta N$ and the mean time between false alarms as $\mathbf{E}_\theta N$. The short ARL can be defined by the following worst case detection delay, proposed by Lorden [14]:

$$\overline{\mathbf{E}}_\lambda N = \sup_{\nu \geq 1}(\operatorname{ess\,sup} \mathbf{E}_{\theta,\lambda}^{(\nu)}[(N - \nu + 1)^+ | X_1, \ldots, X_{\nu-1}]).$$

Note that the definition of $\overline{\mathbf{E}}_\lambda N$ does not depend upon the pre-change distribution $f_\theta$ by virtue of the essential supremum, which takes the "worst possible $X$'s before the change." In our theorems we can also use the average detection delay, proposed by Shiryayev [25] and Pollak [19], $\sup_{\theta \in \Theta}(\sup_{\nu \geq 1} \mathbf{E}_{\theta,\lambda}^\nu(N - \nu | N \geq \nu))$, which is asymptotically equivalent to $\overline{\mathbf{E}}_\lambda N$.

If $\Theta$ and $\Lambda$ are simple, say $\Theta = \{\theta\}$ and $\Lambda = \{\lambda\}$, Page's CUSUM procedure is defined by

$$(1.1) \qquad T_{\mathrm{CM}}(\theta, a) = \inf\left\{n \geq 1 : \max_{1 \leq k \leq n} \sum_{i=k}^n \log \frac{g_\lambda(X_i)}{f_\theta(X_i)} \geq a\right\},$$

where the notation is used to emphasize that the pre-change distribution is $f_\theta$. Moustakides [16] and Ritov [23] showed that Page's CUSUM procedure $T_{\mathrm{CM}}(\theta, a)$ is exactly optimal in the following minimax sense: For any $a > 0$, $T_{\mathrm{CM}}(\theta, a)$ minimizes $\overline{\mathbf{E}}_\lambda N$ among all stopping times $N$ satisfying $\mathbf{E}_\theta N \geq$ $\mathbf{E}_\theta T_{\mathrm{CM}}(\theta, a)$. Earlier Lorden [14] proved this property holds asymptotically. Specifically, Lorden [14] showed that for each pair $(\theta, \lambda)$

$$(1.2) \qquad \overline{\mathbf{E}}_\lambda N \geq (1 + o(1)) \frac{\log \mathbf{E}_\theta N}{I(\lambda, \theta)},$$

as $\mathbf{E}_\theta N \to \infty$ and $T_{\mathrm{CM}}(\theta, a)$ attains the lower bound asymptotically. Here $I(\lambda, \theta) = \mathbf{E}_\lambda \log(g_\lambda(X)/f_\theta(X))$ is the Kullback–Leibler information number. This suggests defining the asymptotic efficiency of a family $\{N(a)\}$ as

$$(1.3) \qquad e(\theta, \lambda) = \liminf_{a \to \infty} \frac{\log \mathbf{E}_\theta N(a)}{I(\lambda, \theta) \overline{\mathbf{E}}_\lambda N(a)},$$

where $\{N(a)\}$ is required to satisfy $\mathbf{E}_\theta N(a) \to \infty$ as $a \to \infty$. Then $e(\theta, \lambda) \leq 1$ for all families, so we can define:

DEFINITION 1.1. A family of stopping times $\{N(a)\}$ is *asymptotically efficient* at $(\theta, \lambda)$ if $e(\theta, \lambda) = 1$.

It follows that Page's CUSUM procedure $T_{\mathrm{CM}}(\theta, a)$ for detecting a change in distribution from $f_\theta$ to $g_\lambda$ is asymptotically efficient at $(\theta, \lambda)$. However, $T_{\mathrm{CM}}(\theta, a)$ in general will not be asymptotically efficient at $(\theta', \lambda)$ if $\theta' \neq \theta$; see Section 2.4 in [31], equation (2.57) in [28] and Table 1 in [5].

**2. Simple post-change hypotheses.** It will be assumed in this section and only in this section that $f_\theta$ and $g_\lambda = f_\lambda$ belong to a one-parameter exponential family

$$(2.1) \qquad f_\xi(x) = \exp(\xi x - b(\xi)), \qquad -\infty < x < \infty, \xi \in \Omega,$$

with natural parameter space $\Omega = (\underline{\xi}, \bar{\xi})$ with respect to a $\sigma$-finite measure $F$. Then $b(\xi)$ is strictly convex on $\Omega$. Assume that $\Theta = [\theta_0, \theta_1]$ is a subset of $\Omega$, and $\lambda$ is a given value outside the interval $\Theta$, say $\lambda > \theta_1$. In this section we consider the problem of detecting a change in distribution from $f_\theta$ for some $\theta \in \Theta$ to $f_\lambda$ and we want to find a stopping time $N$ such that $\mathbf{E}_\theta N$ is as large as possible for each $\theta \in \Theta = [\theta_0, \theta_1]$ subject to the constraint

$$(2.2) \qquad \overline{\mathbf{E}}_\lambda N \leq \gamma,$$

where $\gamma > 0$ is a given constant and $\lambda \notin \Theta$.

One cannot simultaneously maximize $\mathbf{E}_\theta N$ for all $\theta \in \Theta$ subject to (2.2) since the maximum for each $\theta$ is uniquely attained by Page's CUSUM procedure $T_{\mathrm{CM}}(\theta, a)$ in (1.1). As one referee pointed out, if one wants to maximize $\inf_{\theta \in \Theta} \mathbf{E}_\theta N$ subject to (2.2), then the exactly optimal solution is Page's CUSUM procedure $T_{\mathrm{CM}}(\theta_1, a)$ for detecting a change in distribution from $f_{\theta_1}$ to $f_\lambda$. This is because $\inf_{\theta \in \Theta} \mathbf{E}_\theta N \leq \mathbf{E}_{\theta_1} N$ with equality holding for $N = T_{\text{CM}}(\theta_1, a)$, which maximizes $\mathbf{E}_{\theta_1} N$ among all stopping times $N$ satisfying $\bar{\mathbf{E}}_\lambda N \leq \bar{\mathbf{E}}_\lambda T_{\text{CM}}(\theta_1, a)$. In other words, this setup is equivalent to the simplest problem of detecting a change in distribution from $f_{\theta_1}$ to $f_\lambda$.

In this section, rather than be satisfied with just $\inf_{\theta \in \Theta} \mathbf{E}_\theta N$, a lower bound on $\mathbf{E}_\theta N$ over $\theta \in \Theta$, we want to maximize $\mathbf{E}_\theta N$ asymptotically for each $\theta \in \Theta$ as $\gamma \to \infty$, or equivalently, to find a family of stopping times that is asymptotically efficient at $(\theta, \lambda)$ for every $\theta \in \Theta = [\theta_0, \theta_1]$.

Before studying change-point problems in Section 2.2, we first consider the corresponding open-ended hypothesis testing problems in Section 2.1, since the basic arguments are clearer for hypothesis testing problems and are readily extendable to change-point problems.

2.1. *Open-ended hypothesis testing.* Suppose $X_1, X_2, \ldots$ are independent and identically distributed random variables with probability density $f_\xi$ of the form (2.1) on the natural parameter space $\Omega = (\underline{\xi}, \bar{\xi})$. Suppose we are interested in testing the null hypothesis

$$H_0 : \xi \in \Theta = [\theta_0, \theta_1]$$

against the alternative hypothesis

$$H_1 : \xi \in \Lambda = \{\lambda\},$$

where $\underline{\xi} < \theta_0 < \theta_1 < \lambda < \bar{\xi}$.

Motivated by applications to change-point problems, we consider the following open-ended hypothesis testing problems. Assume that if $H_0$ is true, sampling costs nothing and our preferred action is just to observe $X_1, X_2, \ldots$ without stopping. On the other hand, if $H_1$ is true, each observation costs a fixed amount and we want to stop sampling as soon as possible and reject the null hypothesis $H_0$.

Since there is only one terminal decision, a statistical procedure for an open-ended hypothesis testing problem is defined by a stopping time $N$. The null hypothesis $H_0$ is rejected if and only if $N < \infty$. A good procedure $N$ should keep the error probabilities $\mathbf{P}_\theta(N < \infty)$ small for every $\theta \in \Theta$ while keeping $\mathbf{E}_\lambda N$ small.

The problem in this subsection is to find a stopping time $N$ such that $\mathbf{P}_\theta(N < \infty)$ will be as small as possible for every $\theta \in \Theta = [\theta_0, \theta_1]$ subject to the constraint

$$(2.3) \qquad\qquad\qquad \mathbf{E}_\lambda N \leq \gamma,$$

where $\gamma > 0$ is a given constant.

For each $\theta \in \Theta$, by [32], the minimum of $\mathbf{P}_\theta(N < \infty)$ is uniquely attained by the one-sided sequential probability ratio test (SPRT) of $H_{0,\theta} : \xi = \theta$ versus $H_1 : \xi = \lambda$, which is given by

$$\tau_\theta = \inf\left\{ n \geq 1 : \sum_{i=1}^n \log \frac{f_\lambda(X_i)}{f_\theta(X_i)} \geq C_\theta \right\}.$$

In order to satisfy (2.3), it is well known that $C_\theta \approx I(\lambda, \theta)\gamma$; see, for example, page 26 of [28]. A simple observation is that the null hypothesis is expressed as a union of the individual null hypotheses, $H_{0,\theta}\colon \xi = \theta$, and so the intersection-union method (see [2]) suggests considering the stopping time

$$(2.4) \quad M(a) = \inf\left\{ n \geq 1 : \sum_{i=1}^{n} \log \frac{f_\lambda(X_i)}{f_\theta(X_i)} \geq I(\lambda, \theta)a \text{ for all } \theta_0 \leq \theta \leq \theta_1 \right\}.$$

The rationale is that $H_0$ can be rejected only if each of the individual null hypotheses $H_{0,\theta}\colon \xi = \theta$ can be rejected.

In order to study the behavior of $M(a)$, it is useful to express $M(a)$ in terms of $S_n = X_1 + \cdots + X_n$. Define

$$(2.5) \quad \phi(\theta) = \frac{b(\lambda) - b(\theta)}{\lambda - \theta}.$$

Then by (2.1), the stopping time $M(a)$ can be written as

$$(2.6) \quad M(a) = \inf\left\{ n \geq 1 : S_n \geq b'(\lambda)a + \sup_{\theta_0 \leq \theta \leq \theta_1}[(n-a)\phi(\theta)] \right\}$$

because $\lambda > \theta_1$. Now $\phi(\theta)$ is an increasing function since $b(\theta)$ is convex, thus the supremum in (2.6) is attained at $\theta = \theta_0$ if $n \leq a$, and at $\theta = \theta_1$ if $n > a$. Therefore, $M(a)$ is equivalent to the simpler test which uses two simultaneous SPRTs (with appropriate boundaries), one for each of the individual null hypotheses $\theta_0, \theta_1$. This fact makes it convenient for theoretical analysis and numerical simulations.

The following theorem, whose proof is given in Section 5, establishes the asymptotic properties of $M(a)$ for large $a$.

THEOREM 2.1. *For any $a > 0$ and all $\theta_0 \leq \theta \leq \theta_1$*

$$(2.7) \quad \frac{|\log \mathbf{P}_\theta(M(a) < \infty)|}{I(\lambda, \theta)} \geq a,$$

*and as $a \to \infty$*

$$(2.8) \quad \mathbf{E}_\lambda M(a) = a + (C + o(1))\sqrt{a},$$

*where*

$$(2.9) \quad C = \left( \frac{\lambda - \theta_1}{I(\lambda, \theta_1)} - \frac{\lambda - \theta_0}{I(\lambda, \theta_0)} \right)\sqrt{\frac{b''(\lambda)}{2\pi}} > 0.$$

The following corollary establishes the asymptotic optimality of $M(a)$.

COROLLARY 2.1. *Suppose* $\{N(a)\}$ *is a family of stopping times such that* $\mathbf{E}_\lambda N(a) \le \mathbf{E}_\lambda M(a)$. *For all* $\theta_0 \le \theta \le \theta_1$ *as* $a \to \infty$,

$$\frac{|\log \mathbf{P}_\theta(N(a) < \infty)|}{I(\lambda, \theta)} \le a + (C + o(1))\sqrt{a},$$

*where* $C$ *is as defined in* (2.9). *Thus* $M(a)$ *asymptotically minimizes the error probabilities* $\mathbf{P}_\theta(N < \infty)$ *for every* $\theta \in \Theta = [\theta_0, \theta_1]$ *among all stopping times* $N$ *such that* $\mathbf{E}_\lambda N \le \mathbf{E}_\lambda M(a)$.

PROOF. The corollary follows directly from Theorem 2.1 and the well-known fact that

$$\frac{|\log \mathbf{P}_\theta(N(a) < \infty)|}{I(\lambda, \theta)} \le \mathbf{E}_\lambda N(a)$$

for all $\theta \in [\theta_0, \theta_1]$. □

2.2. *Change-point problems.* Now let us consider the problem of detecting a change in distribution from $f_\theta$ for some $\theta \in \Theta = [\theta_0, \theta_1]$ to $f_\lambda$. As described earlier, we seek a family of stopping times that is asymptotically efficient at $(\theta, \lambda)$ for every $\theta \in \Theta$.

A method for finding such a family is suggested by the following result, which indicates the relationship between open-ended hypothesis testing and change-point problems.

LEMMA 2.1 (Lorden [14]). *Let* $N$ *be a stopping time with respect to* $X_1, X_2, \ldots$. *For* $k = 1, 2, \ldots$, *let* $N_k$ *denote the stopping time obtained by applying* $N$ *to* $X_k, X_{k+1}, \ldots$ *for* $k = 1, 2, \ldots$, *and define*

$$N^* = \min_{k \ge 1}(N_k + k - 1).$$

*Then* $N^*$ *is a stopping time with*

$$\mathbf{E}_\theta N^* \ge 1/\mathbf{P}_\theta(N < \infty) \quad and \quad \overline{\mathbf{E}}_\lambda N^* \le \mathbf{E}_\lambda N$$

*for any* $\theta$ *and* $\lambda$.

Let $M(a)$ be the stopping time defined in (2.4), and let $M_k(a)$ be the stopping time obtained by applying $M(a)$ to the observations $X_k, X_{k+1}, \ldots$. Define a new stopping time by $M^*(a) = \min_{k \ge 1}(M_k(a) + k - 1)$. In other words,

$$(2.10) \quad M^*(a) = \inf\left\{ n \ge 1 : \max_{1 \le k \le n} \inf_{\theta_0 \le \theta \le \theta_1}\left( \sum_{i=k}^n \log \frac{f_\lambda(X_i)}{f_\theta(X_i)} - I(\lambda, \theta)a \right) \ge 0 \right\}.$$

The next theorem establishes the asymptotic performance of $M^*(a)$, which immediately implies that the family $\{M^*(a)\}$ is asymptotically efficient at $(\theta, \lambda)$ for every $\theta \in \Theta$.

THEOREM 2.2. *For any $a > 0$ and $\theta_0 \le \theta \le \theta_1$,*

$$\mathbf{E}_\theta M^*(a) \ge \exp(I(\lambda, \theta)a), \tag{2.11}$$

*and as $a \to \infty$,*

$$\overline{\mathbf{E}}_\lambda M^*(a) \le a + (C + o(1))\sqrt{a}, \tag{2.12}$$

*where $C$ is as defined in (2.9). Moreover, if $\{N(a)\}$ is a family of stopping times such that (2.11) holds for some $\theta$ with $N(a)$ replacing $M^*(a)$, then*

$$\overline{\mathbf{E}}_\lambda N(a) \ge a + O(1) \qquad as \ a \to \infty. \tag{2.13}$$

PROOF. Relations (2.11) and (2.12) follow at once from Theorem 2.1 and Lemma 2.1. Relation (2.13) follows from the following proposition, which improves Lorden's lower bound in (1.2). $\square$

PROPOSITION 2.1. *Given $\theta$ and $\lambda \ne \theta$, there exists an $M = M(\theta, \lambda) > 0$ such that for any stopping time $N$,*

$$\log \mathbf{E}_\theta N \le I(\lambda, \theta)\overline{\mathbf{E}}_\lambda N + M. \tag{2.14}$$

PROOF. By equation (2.53) on page 26 of [28], there exist $C_1$ and $C_2$ such that for Page's CUSUM procedure $T_{\mathrm{CM}}(\theta, a)$ in (1.1),

$$\mathbf{E}_\theta T_{\mathrm{CM}}(\theta, a) \le C_1 e^a \quad \text{and} \quad I(\lambda, \theta)\overline{\mathbf{E}}_\lambda T_{\mathrm{CM}}(\theta, a) \ge a - C_2$$

for all $a > 0$. For any given stopping time $N$, choose $a = \log \mathbf{E}_\theta N - \log C_1$; then $\mathbf{E}_\theta N = C_1 e^a \ge \mathbf{E}_\theta T_{\mathrm{CM}}(\theta, a)$. The optimality property of $T_{\mathrm{CM}}(\theta, a)$ [16] implies that

$$I(\lambda, \theta) \log \overline{\mathbf{E}}_\lambda N \ge I(\lambda, \theta) \log \overline{\mathbf{E}}_\lambda T_{\mathrm{CM}}(\theta, a)$$
$$\ge a - C_2 = \log \mathbf{E}_\theta N - \log C_1 - C_2. \qquad \square$$

The following corollary follows at once from Theorem 2.2.

COROLLARY 2.2. *Suppose $\{N(a)\}$ is a family of stopping times such that*

$$\overline{\mathbf{E}}_\lambda N(a) \le \overline{\mathbf{E}}_\lambda M^*(a).$$

*Then for all $\theta_0 \le \theta \le \theta_1$, as $a \to \infty$,*

$$\frac{\log \mathbf{E}_\theta N(a)}{I(\lambda, \theta)} \le a + (C + o(1))\sqrt{a},$$

*where $C$ is as defined in (2.9). Thus, as $a \to \infty$, $M^*(a)$ asymptotically maximizes $\log \mathbf{E}_\theta N$ [up to $O(\sqrt{a})$] for every $\theta \in [\theta_0, \theta_1]$ among all stopping times $N$ such that $\overline{\mathbf{E}}_\lambda N \le \overline{\mathbf{E}}_\lambda M^*(a)$.*

REMARK. The $O(\sqrt{a})$ terms are the price one must pay for optimality at every pre-change distribution $f_\theta$.

In order to implement the stopping times $M^*(a)$ numerically, using (2.6), we can express $M^*(a)$ in the following convenient form:

$$
M^*(a) = \inf \Bigg\{ n \geq 1 : \max_{n-b+1 \leq k \leq n} \sum_{i=k}^{n} \log \frac{f_\lambda(X_i)}{f_{\theta_0}(X_i)} \geq I(\lambda, \theta_0)a,
$$
(2.15)
$$
\text{or } W_{n-b} + \sum_{i=n-b+1}^{n} \log \frac{f_\lambda(X_i)}{f_{\theta_1}(X_i)} \geq I(\lambda, \theta_1)a \Bigg\},
$$

where $b = [a]$, $W_k = \max\{W_{k-1}, 0\} + \log(f_\lambda(X_k)/f_{\theta_1}(X_k))$ and $W_0 = 0$. Since $W_k$ can be calculated recursively, this form reduces the memory requirements at every stage $n$ from the full data set $\{X_1, \ldots, X_n\}$ to the data set of size $b+1$, that is, $\{X_{n-b}, X_{n-b+1}, \ldots, X_n\}$. It is easy to see that this form involves only $O(a)$ computations at every stage $n$.

As an Associate Editor noted, there are other procedures that can have the same asymptotic optimality properties as $M^*(a)$. For example, if we define a slightly different procedure $M_1^*(a)$ by switching $\inf_{\theta_0 \leq \theta \leq \theta_1}$ with $\max_{1 \leq k \leq n}$ in the definition of $M^*(a)$ in (2.10), or if we define $M_2^*(a) = \sup_{\theta_0 \leq \theta \leq \theta_1}\{T_{\mathrm{CM}}(\theta, I(\lambda, \theta)a)\}$, where $T_{\mathrm{CM}}(\theta, I(\lambda, \theta)a)$ is Page's CUSUM procedure for detecting a change in distribution from $f_\theta$ to $f_\lambda$ with log-likelihood ratio boundary $I(\lambda, \theta)a$, then both $M_1^*(a)$ and $M_2^*(a)$ are well-defined stopping times that are asymptotically efficient at $(\theta, \lambda)$ for every $\theta \in \Theta$. However, both $M_1^*(a)$ and $M_2^*(a)$ are difficult to implement, although one can easily implement their approximations which replace $\Theta = [\theta_0, \theta_1]$ by a (properly chosen) finite subset of $\Theta$.

It is important to emphasize that in all the above procedures we should choose appropriate stopping boundaries. Otherwise the procedures may not be asymptotically efficient at every $\theta \in \Theta$. For instance, motivated by the generalized likelihood ratio method, one may want to use the procedure

$$
T'(a) = \inf \Bigg\{ n \geq 1 : \max_{1 \leq k \leq n} \log \frac{g_\lambda(X_k) \cdots g_\lambda(X_n)}{\sup_{\theta_0 \leq \theta \leq \theta_1}(f_\theta(X_k) \cdots f_\theta(X_n))} \geq a \Bigg\}
$$
$$
= \inf \Bigg\{ n \geq 1 : \max_{1 \leq k \leq n} \inf_{\theta_0 \leq \theta \leq \theta_1} \Bigg[ (\lambda - \theta) \sum_{i=k}^{n}(X_i - \phi(\theta)) \Bigg] \geq a \Bigg\},
$$

where $\phi(\theta)$ is defined in (2.5). Unfortunately, for all $a > 0$, $T'(a)$ is equivalent to Page's CUSUM procedure $T_{\mathrm{CM}}(\theta_1, a)$, and thus it will not be asymptotically efficient at every $\theta$. To see this, first note that $T'(a) \geq T_{\mathrm{CM}}(\theta_1, a)$ by their definitions. Next, if $T_{\mathrm{CM}}(\theta_1, a)$ stops at time $n_0$, then for some $1 \leq k_0 \leq n_0$, $\sum_{i=k_0}^{n_0}(X_i - \phi(\theta_1)) \geq a/(\lambda - \theta_1)$ since $\lambda > \theta_1$. Thus, if $a > 0$, then for all $\theta_0 \leq \theta \leq \theta_1 (\leq \lambda)$,

$$\sum_{i=k_0}^{n_0}(X_i - \phi(\theta)) \geq \sum_{i=k_0}^{n_0}(X_i - \phi(\theta_1)) \geq \frac{a}{\lambda - \theta_1} \geq \frac{a}{\lambda - \theta}$$

because $\phi(\theta)$ is an increasing function of $\theta$. This implies that $T'(a)$ stops before or at time $n_0$ and so $T'(a) \leq T_{\mathrm{CM}}(\theta_1, a)$. Therefore, $T'(a) = T_{\mathrm{CM}}(\theta_1, a)$. Similarly, if one considers $T''(a) = \sup_{\theta_0 \leq \theta \leq \theta_1}\{T_{\mathrm{CM}}(\theta, a)\}$, then $T''(a)$ is also equivalent to $T_{\mathrm{CM}}(\theta_1, a)$, because for all $a > 0$, Page's CUSUM procedure $T_{\mathrm{CM}}(\theta, a)$ is increasing as a function of $\theta \in [\theta_0, \theta_1]$ in the sense that $T_{\mathrm{CM}}(\theta, a) \leq T_{\mathrm{CM}}(\theta', a)$ if $\theta \leq \theta'$.

2.3. *Extension to half-open interval.* Suppose $X_1, X_2, \ldots$ are independent and identically distributed random variables with probability density $f_\xi$ of the form (2.1) and suppose we are interested in testing the null hypothesis

$$H_0 : \xi \in \Theta = (\underline{\xi}, \theta_1]$$

against the alternative hypothesis

$$H_1 : \xi \in \Lambda = \{\lambda\},$$

where $\theta_1 < \lambda$. Recall that $\Omega = (\underline{\xi}, \bar{\xi})$ is the natural parameter space of $\xi$. Assume

$$(2.16) \qquad \lim_{\theta \to \underline{\xi}} \mathbf{E}_\theta X = -\infty.$$

This condition is equivalent to $\lim_{\theta \to \underline{\xi}} b'(\theta) = -\infty$ since $b'(\theta) = \mathbf{E}_\theta X$. Many distributions satisfy this condition. For example, (2.16) holds for the normal distributions since $\mathbf{E}_\theta X = \theta$ and $\underline{\xi} = -\infty$. It also holds for the negative exponential density since $b(\theta) = -\log \theta, \underline{\xi} = 0$ and $\mathbf{E}_\theta X = b'(\theta) = -1/\theta$.

As in (2.4), our proposed open-ended test $M(a)$ of $H_0 : \xi \in \Theta = (\underline{\xi}, \theta_1]$ against $H_1 : \xi = \lambda$ is defined by

$$\hat{M}(a) = \inf\left\{ n \geq 1 : \sum_{i=1}^{n} \log \frac{f_\lambda(X_i)}{f_\theta(X_i)} \geq I(\lambda, \theta)a \text{ for all } \underline{\xi} \leq \theta \leq \theta_1 \right\}.$$

As in (2.6), $\hat{M}(a)$ can be written as

$$(2.17) \quad \hat{M}(a) = \inf\left\{ n \geq 1 : \sum_{i=1}^{n} X_i \geq b'(\lambda)a + \sup_{\underline{\xi} < \theta \leq \theta_1} [(n-a)\phi(\theta)] \right\},$$

where $\phi(\theta)$ is defined in (2.5). By L'Hôpital's rule and the condition in (2.16),

$$\lim_{\theta \to \underline{\xi}} \phi(\theta) = \lim_{\theta \to \underline{\xi}} \frac{b(\lambda) - b(\theta)}{\lambda - \theta} = \lim_{\theta \to \underline{\xi}} b'(\theta) = \lim_{\theta \to \underline{\xi}} \mathbf{E}_\theta X = -\infty.$$

Thus for any $n < a$, $\sum_{i=1}^{n} X_i$ is finite but $\sup_{\underline{\xi} < \theta \le \theta_1}[(n-a)\phi(\theta)] = \infty$. So $\hat{M}(a)$ will never stop at time $n < a$. Recall that $\phi(\theta)$ is an increasing function of $\theta$, hence the supremum in (2.17) is attained at $\theta = \theta_1$ if $n \ge a$. Therefore,

$$\hat{M}(a) = \inf\left\{ n \ge a : \sum_{i=1}^{n} \log \frac{f_\lambda(X_i)}{f_{\theta_1}(X_i)} \ge I(\lambda, \theta_1)a \right\}.$$

For the problem of detecting a change in distribution from some $f_\theta$ with $\theta \in \Theta = (\underline{\xi}, \theta_1]$ to $f_\lambda$, define $\hat{M}^*(a)$ from $\hat{M}(a)$ as before, so that

$$\hat{M}^*(a) = \inf\left\{ n \ge a : \max_{1 \le k \le n-a+1} \sum_{i=k}^{n} \log \frac{f_\lambda(X_i)}{f_{\theta_1}(X_i)} \ge I(\lambda, \theta_1)a \right\}.$$

Using arguments similar to the proof of Theorem 2.2, we have:

THEOREM 2.3. *For $a > 0$ and $\theta \in (\underline{\xi}, \theta_1]$,*

$$\mathbf{E}_\theta \hat{M}^*(a) \ge \exp(I(\lambda, \theta)a),$$

*and as $a \to \infty$,*

$$\overline{\mathbf{E}}_\lambda \hat{M}^*(a) \le a + (\hat{C} + o(1))\sqrt{a},$$

*where*

$$\hat{C} = \frac{\lambda - \theta_1}{I(\lambda, \theta_1)} \sqrt{\frac{b''(\lambda)}{2\pi}} > 0.$$

Thus the analogue of Corollary 2.2 holds, and so $\hat{M}^*(a)$ asymptotically maximizes $\log \mathbf{E}_\theta N$ [up to $O(\sqrt{a})$] for every $\theta \in (\underline{\xi}, \theta_1]$ among all stopping times $N$ such that $\overline{\mathbf{E}}_\lambda N \le \overline{\mathbf{E}}_\lambda \hat{M}^*(a)$.

TABLE 1
*Long ARL for different procedures*

| $\theta$ | Best possible | $M^*(a)$ $(a = 18.50)$ | $T_{\mathrm{CM}}(-0.5, a)$ $(a = 2.92)$ | $T_{\mathrm{CM}}(-1.0, a)$ $(a = 9.88)$ |
|---|---|---|---|---|
| $-0.5$ | $233 \pm 7$ | $206 \pm 6$ | $233 \pm 7$ | $125 \pm 3$ |
| $-0.6$ | $523 \pm 15$ | $501 \pm 15$ | $518 \pm 15$ | $297 \pm 8$ |
| $-0.7$ | $1384 \pm 43$ | $1324 \pm 43$ | $1227 \pm 37$ | $938 \pm 29$ |
| $-0.8$ | $5157 \pm 165$ | $4688 \pm 148$ | $3580 \pm 113$ | $4148 \pm 129$ |
| $-0.9$ | $22{,}942 \pm 699$ | $19{,}217 \pm 606$ | $10{,}613 \pm 343$ | $21{,}617 \pm 658$ |
| $-1.0$ | $118{,}223 \pm 3711$ | $83{,}619 \pm 2566$ | $31{,}641 \pm 1036$ | $118{,}223 \pm 3711$ |

(The best possible values are obtained from an optimal envelope of Page's CUSUM procedures.)

2.4. *Numerical examples.* In this subsection we describe the results of a Monte Carlo experiment designed to check the insights obtained from the asymptotic theory of previous subsections. The simulations consider the problem of detecting a change in a normal mean, where the pre-change distribution $f_\theta = N(\theta, 1)$ with $\theta \in \Theta = [-1, -0.5]$, and the post-change distribution $f_\lambda = N(\lambda, 1)$ with $\lambda \in \Lambda = \{0\}$.

Table 1 compares our procedure $M^*(a)$ and two versions of Page's CUSUM procedure $T_{\mathrm{CM}}(\theta_0, a)$ over a range of $\theta$ values. Here

$$T_{\mathrm{CM}}(\theta_0, a) = \inf\left\{n \geq 1 : \max_{1 \leq k \leq n} \sum_{i=k}^n \log \frac{f_\lambda(X_i)}{f_{\theta_0}(X_i)} \geq a\right\}$$

$$= \inf\left\{n \geq 1 : \max_{1 \leq k \leq n} \sum_{i=k}^n (-\theta_0)\left[X_i - \frac{\theta_0}{2}\right] \geq a\right\}.$$

The threshold value $a$ for Page's CUSUM procedure $T_{\mathrm{CM}}(\theta_0, a)$ and our procedure $M^*(a)$ was determined from the criterion $\overline{\mathbf{E}}_\lambda N \approx 20$. First, a $10^4$-repetition Monte Carlo simulation was performed to determine the appropriate values of $a$ to yield the desired detection delay to within the range of sampling error. With the thresholds used, the detection delay $\overline{\mathbf{E}}_\lambda N$ is close enough to 20 so that the difference is negligible, that is, correcting the threshold to get exactly 20 (if we knew how to do that) would change $\mathbf{E}_\theta N$ by an amount that would make little difference in light of the simulation errors $\mathbf{E}_\theta N$ already has. Next, using the obtained threshold value $a$, we ran 1000 repetitions to simulate long ARL, $\mathbf{E}_\theta N$, for different $\theta$.

Table 1 also reports the best possible $\mathbf{E}_\theta N$ at each of the values of $\theta$ subject to $\overline{\mathbf{E}}_\lambda N \approx 20$. Note that they are obtained from an optimal envelope of Page's CUSUM procedures and therefore cannot be attained simultaneously in practice. Each result in Table 1 is recorded as the Monte Carlo estimate ± standard error.

Table 1 shows that $M^*(a)$ performs well over a broad range of $\theta$, which is consistent with the asymptotic theory of $M^*(a)$ developed in Sections 2.2 and 2.3 showing that $M^*(a)$ attains [up to $O(\sqrt{a})$] the asymptotic upper bounds for $\log \mathbf{E}_\theta N$ in Corollary 2.2 as $a \to \infty$.

**3. Composite post-change hypotheses.** Let $\Theta$ and $\Lambda$ be two compact disjoint subsets of some Euclidean space. Let $\{f_\theta; \theta \in \Theta\}$ and $\{g_\lambda; \lambda \in \Lambda\}$ be two sets of densities, absolutely continuous with respect to the same nondegenerate $\sigma$-finite measure. In this section we are interested in detecting a change in distribution from $f_\theta$ for some $\theta \in \Theta$ to $g_\lambda$ for some $\lambda \in \Lambda$. Here we no longer assume the densities belong to exponential families, and we assume that both $\Theta$ and $\Lambda$ are composite.

Ideally we would like a stopping time $N$ which minimizes the detection delay $\overline{\mathbf{E}}_\lambda N$ for all $\lambda \in \Lambda$ and maximizes $\mathbf{E}_\theta N$ for all $\theta \in \Theta$, that is, we seek a family $\{N(a)\}$ which is asymptotically efficient for all $(\theta, \lambda) \in \Theta \times \Lambda$. However, in general such a family does not exist. For example, for $\Lambda = \{\lambda_1, \lambda_2\}$ it is easy to see from (1.3) that there exists a family that is asymptotically efficient at both $(\theta, \lambda_1)$ and $(\theta, \lambda_2)$ for all $\theta \in \Theta$ only if $I(\lambda_2, \theta)/I(\lambda_1, \theta)$ is constant in $\theta \in \Theta$. This fails in general when $\Theta$ is composite. For example, if $f_\theta$ and $g_\lambda$ belong to a one-parameter exponential family and $\Theta$ is an interval, a simple argument shows that $I(\lambda_2, \theta)/I(\lambda_1, \theta)$ is a constant if and only if $\lambda_1 = \lambda_2$.

It is natural to consider the following definition:

DEFINITION 3.1. A family of stopping times $\{N(a)\}$ is asymptotically optimal to first order if:

(i) for each $\theta \in \Theta$, there exists at least one $\lambda_\theta \in \Lambda$ such that the family is asymptotically efficient at $(\theta, \lambda_\theta)$; and

(ii) for each $\lambda \in \Lambda$, there exists at least one $\theta_\lambda \in \Theta$ such that the family is asymptotically efficient at $(\theta_\lambda, \lambda)$.

REMARK. An equivalent definition is to require that the family $\{N(a)\}$ is asymptotically efficient at $(h_1(\delta), h_2(\delta))$ for $\delta \in \Delta$, where $\theta = h_1(\delta)$ and $\lambda = h_2(\delta)$ are onto (not necessary one-to-one) functions from $\Delta$ to $\Theta$ and $\Lambda$, respectively. It is obvious that the standard formulation with simple $\Theta$ and our formulation in Section 2 are two special cases of this definition.

REMARK. It is worth noting that a family of stopping times that is asymptotically optimal to first order is asymptotically admissible in the following sense. A family of stopping times $\{N(a)\}$ is asymptotically inadmissible if there exists another family of stopping times $\{N'(a)\}$ such that for all $\theta \in \Theta$ and all $\lambda \in \Lambda$,

$$\limsup_{a \to \infty} \frac{\log \mathbf{E}_\theta N(a)}{\log \mathbf{E}_\theta N'(a)} \leq 1 \quad \text{and} \quad \liminf_{a \to \infty} \frac{\overline{\mathbf{E}}_\lambda N(a)}{\overline{\mathbf{E}}_\lambda N'(a)} \geq 1,$$

with strict inequality holding for some $\theta$ or $\lambda$. A family of stopping times is asymptotically admissible if it is not asymptotically inadmissible.

Note that when $\Lambda = \{\lambda\}$ is simple, the asymptotically optimal procedure developed in Section 2 satisfies

$$(3.1) \qquad\qquad \log \mathbf{E}_\theta N(a) \sim I(\lambda, \theta)a \qquad \text{as } a \to \infty.$$

Here and everywhere below, $x(a) \sim y(a)$ as $a \to \infty$ means that $\lim_{a\to\infty}(x(a)/y(a)) = 1$. However, when one considers multiple values of the post-change parameter $\lambda$ it is no longer possible to find a procedure such that (3.1) holds for all $(\theta, \lambda) \in \Theta \times \Lambda$. A natural idea is then to seek procedures such that

$$\log \mathbf{E}_\theta N(a) \sim p(\theta)a,$$

where $p(\theta)$ is suitably chosen. It turns out that for "good" choices of $p(\theta)$ one can define $\{N(a)\}$ to be asymptotically optimal to first order.

To accomplish this, first consider the following definitions.

DEFINITION 3.2.   A positive continuous function $p(\cdot)$ on $\Theta$ is an optimizer if for some positive continuous $q(\cdot)$ on $\Lambda$

$$p(\theta) = \inf_{\lambda \in \Lambda} \frac{I(\lambda, \theta)}{q(\lambda)}.$$

Similarly, $q(\cdot)$ on $\Lambda$ is an optimizer if for some positive continuous $p(\cdot)$ on $\Theta$

$$q(\lambda) = \inf_{\theta \in \Theta} \frac{I(\lambda, \theta)}{p(\theta)}.$$

DEFINITION 3.3.   Positive continuous functions $p(\cdot), q(\cdot)$ on $\Theta, \Lambda$, respectively, are an optimizer pair if for all $\theta$ and $\lambda$

$$(3.2) \qquad p(\theta) = \inf_{\lambda \in \Lambda} \frac{I(\lambda, \theta)}{q(\lambda)} \quad \text{and} \quad q(\lambda) = \inf_{\theta \in \Theta} \frac{I(\lambda, \theta)}{p(\theta)}.$$

The following proposition characterizes the relation between these two definitions.

PROPOSITION 3.1.   *If $(p, q)$ is an optimizer pair, then $p$ and $q$ are optimizers. Conversely, for every optimizer $p$, there is a $q$ such that $(p, q)$ is an optimizer pair, namely,*

$$q(\lambda) = \inf_{\theta \in \Theta} \frac{I(\lambda, \theta)}{p(\theta)}$$

*and, similarly, for every optimizer $q$ one can obtain an optimizer pair $(p, q)$ by defining*

$$p(\theta) = \inf_{\lambda \in \Lambda} \frac{I(\lambda, \theta)}{q(\lambda)}.$$

PROOF.   It is obvious that $p$ and $q$ are optimizers if $(p, q)$ is an optimizer pair. Since everything is symmetric in the roles of $p$ and $q$, we only need to prove that the first equation of (3.2) holds for the case where $q$ is defined after $p$. Now fix $\theta_0 \in \Theta$. On the one hand, since $q(\lambda)$ is defined as the infimum over $\Theta$, we have $q(\lambda) \leq I(\lambda, \theta_0)/p(\theta_0)$, so $p(\theta_0) \leq I(\lambda, \theta_0)/q(\lambda)$ for all $\lambda \in \Lambda$. Thus

$$(3.3) \qquad p(\theta_0) \leq \inf_{\lambda \in \Lambda} \frac{I(\lambda, \theta_0)}{q(\lambda)}.$$

On the other hand, since $p$ is an optimizer by assumption, there exists a function $q_0(\cdot)$ on $\Lambda$ such that

$$p(\theta) = \inf_{\lambda \in \Lambda} \frac{I(\lambda, \theta)}{q_0(\lambda)}.$$

For any $\lambda_0 \in \Lambda$, we have $p(\theta) \leq I(\lambda_0, \theta)/q_0(\lambda_0)$ and so $I(\lambda_0, \theta)/p(\theta) \geq q_0(\lambda_0)$ for all $\theta \in \Theta$. Hence

$$\inf_{\theta \in \Theta} \frac{I(\lambda_0, \theta)}{p(\theta)} \geq q_0(\lambda_0).$$

Observe that the left-hand side is just our definition for $q(\lambda_0)$, and so $q(\lambda_0) \geq q_0(\lambda_0)$. Since $\lambda_0$ is arbitrary, we have $q(\lambda) \geq q_0(\lambda)$ for all $\lambda \in \Lambda$. Thus,

$$\inf_{\lambda \in \Lambda} \frac{I(\lambda, \theta_0)}{q(\lambda)} \leq \inf_{\lambda \in \Lambda} \frac{I(\lambda, \theta_0)}{q_0(\lambda)} = p(\theta_0)$$

by using the definition of $p(\theta)$. The first equation of (3.2) follows at once from this and (3.3). □

In fact, Proposition 3.1 provides a method to construct optimizer pairs. One can start with any positive continuous function $q_0(\lambda)$, get an optimizer $p(\theta)$ from it by (3.2) and use the other part of (3.2) to get a $(p, q)$ optimizer pair. Similarly, one can also get a $(p, q)$ optimizer pair by starting with a $p_0(\theta)$.

Now we can define our proposed procedures based on an optimizer $p(\theta)$. First, let $\eta$ be an a priori distribution fully supported on $\Lambda$. Define an open-ended test $T(a)$ by

$$(3.4) \quad T(a) = \inf \left\{ n : \inf_{\theta \in \Theta} \left[ \frac{1}{p(\theta)} \log \frac{\int_\Lambda [g_\lambda(X_1) \cdots g_\lambda(X_n)] \eta(d\lambda)}{f_\theta(X_1) \cdots f_\theta(X_n)} \right] \geq a \right\}.$$

Then our proposed procedure is defined by $T^*(a) = \min_{k \geq 1}(T_k(a) + k - 1)$, where $T_k(a)$ is obtained by applying $T(a)$ to $X_k, X_{k+1}, \ldots$. Equivalently,

$$(3.5) \quad T^*(a) = \inf \Big\{ n \geq 1 : $$
$$\max_{1 \leq k \leq n} \inf_{\theta \in \Theta} \left[ \frac{1}{p(\theta)} \log \frac{\int_\Lambda [g_\lambda(X_k) \cdots g_\lambda(X_n)] \eta(d\lambda)}{f_\theta(X_k) \cdots f_\theta(X_n)} \right] \geq a \Big\}.$$

We also define a slightly different procedure $T_1^*(a)$ by switching $\inf_{\theta \in \Theta}$ with $\max_{1 \leq k \leq n}$ in the definition of $T^*(a)$.

Our main results in this section are stated in the next theorem and its corollary, which establish the asymptotic optimality properties of $T^*(a)$ and $T_1^*(a)$. The proofs are given in Section 3.1.

THEOREM 3.1. *Assume that Assumptions* A1 *and* A2 *below hold and* $\Theta$ *and* $\Lambda$ *are compact. If* $p(\theta)$ *is an optimizer, then* $\{T^*(a)\}$ *and* $\{T_1^*(a)\}$ *are asymptotically optimal to first order.*

COROLLARY 3.1. *Under the assumptions of Theorem* 3.1, *if* $\{N(a)\}$ *is a family of procedures such that*

$$\limsup_{a\to\infty} \frac{\overline{\mathbf{E}}_\lambda N(a)}{\overline{\mathbf{E}}_\lambda T^*(a)} \le 1 \qquad \text{for all } \lambda \in \Lambda,$$

*then*

$$\limsup_{a\to\infty} \frac{\log \mathbf{E}_\theta N(a)}{\log \mathbf{E}_\theta T^*(a)} \le 1 \qquad \text{for all } \theta \in \Theta.$$

*Similarly, if*

$$\liminf_{a\to\infty} \frac{\log \mathbf{E}_\theta N(a)}{\log \mathbf{E}_\theta T^*(a)} \ge 1 \qquad \text{for all } \theta \in \Theta,$$

*then*

$$\liminf_{a\to\infty} \frac{\overline{\mathbf{E}}_\lambda N(a)}{\overline{\mathbf{E}}_\lambda T^*(a)} \ge 1 \qquad \text{for all } \lambda \in \Lambda.$$

*The same assertions are true if* $T^*(a)$ *is replaced by* $T_1^*(a)$.

REMARK. Corollary 3.1 shows that our procedures $T^*(a)$ and $T_1^*(a)$ are also asymptotically optimal in the following sense: If a family of procedures $\{N(a)\}$ performs asymptotically as well as our procedures (or better) uniformly over $\Theta$, then our procedures perform asymptotically as well as $\{N(a)\}$ (or better) uniformly over $\Lambda$, and the same is true if the roles of $\Theta$ and $\Lambda$ are reversed.

REMARK. Theorem 3.1 and Corollary 3.1 show another asymptotic optimality property of our procedures $T^*(a)$ and $T_1^*(a)$: If the optimizer $p(\theta)$ is constructed from $q_0(\lambda)$ by the first equation of (3.2), then our procedures asymptotically maximize $\mathbf{E}_\theta N$ for every $\theta \in \Theta$ among all stopping times $N$ satisfying

$$q_0(\lambda)\overline{\mathbf{E}}_\lambda N \le \gamma \qquad \text{for all } \lambda \in \Lambda,$$

where $\gamma > 0$ is given. Here $q_0(\lambda) > 0$ can be thought of as the cost per observation of delay if the post-change observations have distribution $g_\lambda$.

REMARK. Instead of $T(a)$ in (3.4), we can also define the following stopping time in open-ended hypothesis testing problems:

$$(3.6) \quad \hat{T}(a) = \inf\left\{ n \ge 1 : \inf_{\theta \in \Theta}\left[ \frac{1}{p(\theta)} \log \frac{\sup_\lambda [g_\lambda(X_1)\cdots g_\lambda(X_n)]}{f_\theta(X_1)\cdots f_\theta(X_n)} \right] \ge a \right\},$$

and then use it to construct the corresponding procedures in change-point problems. When $f_\theta$ and $g_\lambda$ are from the same one-parameter exponential family, we can obtain an upper bound on $\mathbf{P}_\theta(\hat{T}(a) < \infty)$ by equation (13) on page 636 in [15], and so we get a lower bound on the long ARL. The upper bound on detection delay follows from the fact that $\hat{T}(a) \leq T(a)$. These procedures are, therefore, also asymptotically optimal to first order if $f_\theta$ and $g_\lambda$ belong to one-parameter exponential families.

REMARK. Note that if $p(\theta) \equiv 1$, then all of our procedures are just based on generalized likelihood ratios. However, in the case where $p(\theta) \equiv 1$ is not an optimizer, generalized likelihood ratio procedures may not be asymptotically optimal to first order. In fact, they are asymptotically inadmissible since they are dominated by our procedures based on an optimizer $p(\theta)$ which is obtained by starting with $p_0(\theta) \equiv 1$.

Throughout this section we impose the following assumptions on the densities $f_\theta$ and $g_\lambda$.

ASSUMPTION A1. The Kullback–Leibler information numbers $I(\lambda, \theta) = \mathbf{E}_\lambda \log(g_\lambda(X)/f_\theta(X))$ are finite. Furthermore:

(a) $I_0 = \inf_\lambda \inf_\theta I(\lambda, \theta) > 0$,
(b) $I(\lambda, \theta)$ and $I(\lambda) = \inf_\theta I(\lambda, \theta)$ are both continuous in $\lambda$.

ASSUMPTION A2. For all $\theta$, $\lambda$:

(a) $\mathbf{E}_\lambda[\log(g_\lambda(X)/f_\theta(X))]^2 < \infty$,
(b) $\lim_{\rho \to 0} \mathbf{E}_\lambda[\log \sup_{|\theta' - \theta| \leq \rho} f_{\theta'}(X) - \log f_\theta(X)]^2 = 0$,
(c) $\lim_{\lambda' \to \lambda} \mathbf{E}_\lambda[\log g_{\lambda'}(X) - \log g_\lambda(X)]^2 = 0$.

Assumptions A1 and A2 are part of the Assumptions 2 and 3 in [7]. Assumption A1(a) guarantees that $\Theta$ and $\Lambda$ are "separated."

3.1. *Proof of main results.* First we establish the lower bound on the long ARLs of our procedures $T^*(a)$ and $T_1^*(a)$ for any arbitrary positive function $p(\theta)$.

LEMMA 3.1. *For all $a > 0$ and $\theta \in \Theta$,*
$$\log \mathbf{E}_\theta T^*(a) \geq \log \mathbf{E}_\theta T_1^*(a) \geq p(\theta)a.$$

PROOF. Define
$$t(\theta, a) = \inf\left\{ n \geq 1 : \frac{1}{p(\theta)} \log \frac{\int_\Lambda [g_\lambda(X_1) \cdots g_\lambda(X_n)] \eta(d\lambda)}{f_\theta(X_1) \cdots f_\theta(X_n)} \geq a \right\}$$

and $t^*(\theta, a) = \min_{k \geq 1}(t_k(\theta, a) + k - 1)$, where $t_k(\theta, a)$ is obtained by applying $t(\theta, a)$ to $X_k, X_{k+1}, \ldots$. Then it is clear that $T^*(a) \geq T_1^*(a) \geq t^*(\theta, a)$, and hence

$$\mathbf{E}_\theta T^*(a) \geq \mathbf{E}_\theta T_1^*(a) \geq \mathbf{E}_\theta[t^*(\theta, a)].$$

Using Lemma 2.1 and Wald's likelihood ratio identity, we have

$$\mathbf{E}_\theta[t^*(\theta, a)] \geq \frac{1}{\mathbf{P}_\theta(t(\theta, a) < \infty)} \geq \exp(p(\theta)a),$$

which proves the lemma. $\square$

Next we derive an upper bound on the detection delays of our procedures $T^*(a)$ and $T_1^*(a)$.

LEMMA 3.2. *Suppose that Assumptions* A1 *and* A2 *hold and $\Theta$ is compact. If $p(\theta)$ is a positive continuous function (not necessarily an optimizer) on $\Theta$, then for all $\lambda \in \Lambda$,*

$$\overline{\mathbf{E}}_\lambda T_1^*(a) \leq \overline{\mathbf{E}}_\lambda T^*(a) \leq (1 + o(1)) \frac{a}{q(\lambda)}$$

*as $a \to \infty$, where $q(\lambda)$ is defined by*

$$q(\lambda) = \inf_{\theta \in \Theta} \frac{I(\lambda, \theta)}{p(\theta)}.$$

PROOF. By definition, $\overline{\mathbf{E}}_\lambda T_1^*(a) \leq \overline{\mathbf{E}}_\lambda T^*(a) \leq \mathbf{E}_\lambda T(a)$, where $T(a)$ is defined in (3.4), so it suffices to show that

$$\mathbf{E}_\lambda T(a) \leq (1 + o(1)) \frac{a}{q(\lambda)}$$

for any $\lambda \in \Lambda$. We will use the method in [7] to prove this inequality. Fix $\lambda_0 \in \Lambda$ and choose an arbitrary $\varepsilon > 0$. By Assumptions A1 and A2, the compactness of $\Theta$ and the continuity of $p(\theta)$, there exist a finite covering $\{U_i, 1 \leq i \leq k_\varepsilon\}$ of $\Theta$ (with $\theta_i \in U_i$) and positive numbers $\delta_\varepsilon$ such that for all $\lambda \in V_\varepsilon = \{\lambda \mid |\lambda - \lambda_0| < \delta_\varepsilon\}$, and $i = 1, \ldots, k_\varepsilon$,

$$(3.7) \qquad \mathbf{E}_{\lambda_0}\left[\log g_\lambda(X) - \log \sup_{\theta \in U_i} f_\theta(X)\right] \geq I(\lambda_0, \theta_i) - \varepsilon$$

and

$$\sup_{\theta \in U_i} |p(\theta) - p(\theta_i)| < \varepsilon.$$

Let $N_1(a)$ be the smallest $n$ such that

$$(3.8) \quad \log \int_{V_\varepsilon} [g_\lambda(X_1) \cdots g_\lambda(X_n)]\eta(d\lambda) \geq \sup_{\theta \in \Theta}\left[p(\theta)a + \sum_{j=1}^n \log f_\theta(X_j)\right].$$

Clearly $N_1(a) \geq T(a)$. By Jensen's inequality, the left-hand side of (3.8) is greater than or equal to

$$(3.9) \quad \int_{V_\varepsilon} \log[g_\lambda(X_1) \cdots g_\lambda(X_n)] \frac{\eta(d\lambda)}{\eta(V_\varepsilon)} + \log \eta(V_\varepsilon)$$

$$= \sum_{j=1}^n \int_{V_\varepsilon} \log g_\lambda(X_j) \frac{\eta(d\lambda)}{\eta(V_\varepsilon)} - |\log \eta(V_\varepsilon)|$$

since $\eta(V_\varepsilon) \leq 1$. Since $\{U_i\}$ covers $\Theta$, the right-hand side of (3.8) is less than or equal to

$$\max_{1 \leq i \leq k_\varepsilon} \sup_{\theta \in U_i} \left[ p(\theta)a + \sum_{j=1}^n \log f_\theta(X_j) \right]$$

$$(3.10) \quad \leq \max_{1 \leq i \leq k_\varepsilon} \left[ (p(\theta_i) + \varepsilon)a + \sup_{\theta \in U_i} \sum_{j=1}^n \log f_\theta(X_j) \right]$$

$$\leq \max_{1 \leq i \leq k_\varepsilon} \left[ (p(\theta_i) + \varepsilon)a + \sum_{j=1}^n \log \sup_{\theta \in U_i} f_\theta(X_j) \right].$$

For $j = 1, 2, \ldots$, put

$$Y_j = \int_{V_\varepsilon} \log g_\lambda(X_j) \frac{\eta(d\lambda)}{\eta(V_\varepsilon)} \quad \text{and} \quad Z_j^i = \log \sup_{\theta \in U_i} f_\theta(X_j) \qquad \text{for } i = 1, \ldots, k_\varepsilon.$$

Let $N_2(a)$ be the smallest $n$ such that

$$\sum_{j=1}^n Y_j - \max_{1 \leq i \leq k_\varepsilon} \left[ \sum_{j=1}^n Z_j^i + (p(\theta_i) + \varepsilon)a \right] \geq |\log \eta(V_\varepsilon)|$$

or, equivalently, the smallest $n$ such that for all $1 \leq i \leq k_\varepsilon$,

$$\sum_{j=1}^n \frac{Y_j - Z_j^i}{p(\theta_i)} \geq a \left[ 1 + \frac{\varepsilon}{p(\theta_i)} \right] + \frac{|\log \eta(V_\varepsilon)|}{p(\theta_i)}.$$

Using (3.9) and (3.10), it is clear that $N_2(a) \geq N_1(a)$. Let $p_0 = \inf_{\theta \in \Theta} p(\theta)$; then $p_0 > 0$ since $p(\theta)$ is a positive continuous function and $\Theta$ is compact. Define $\tau_\varepsilon = |\log \eta(V_\varepsilon)|/p_0$, and let $N_3(a)$ be the smallest $n$ such that

$$\min_{1 \leq i \leq k_\varepsilon} \sum_{j=1}^n \frac{Y_j - Z_j^i}{p(\theta_i)} \geq a \left( 1 + \frac{\varepsilon}{p_0} \right) + \tau_\varepsilon$$

or, equivalently,

$$\sum_{j=1}^n \left[ \frac{Y_j - Z_j^1}{p(\theta_1)} - \varepsilon \right] + \min_{1 \leq i \leq k_\varepsilon} \sum_{j=1}^n \left[ \frac{Y_j - Z_j^i}{p(\theta_i)} - \frac{Y_j - Z_j^1}{p(\theta_1)} + \varepsilon \right] \geq a \left( 1 + \frac{\varepsilon}{p_0} \right) + \tau_\varepsilon.$$

Clearly $N_3(a) \geq N_2(a)$. From (3.7) we have

$$(3.11) \qquad \mathbf{E}_{\lambda_0} \left[ \frac{Y_j - Z_j^i}{p(\theta_i)} - \varepsilon \right] \geq \frac{I(\lambda_0, \theta_i)}{p(\theta_i)} - \varepsilon \left( 1 + \frac{1}{p_0} \right) \qquad \text{for } i = 1, \ldots, k_\varepsilon.$$

For $n = 1, 2, \ldots$ define

$$S_n = \sum_{j=1}^{n} \left[ \frac{Y_j - Z_j^1}{p(\theta_1)} - \varepsilon \right]$$

and

$$B_n^i = \sum_{j=1}^{n} \left[ \frac{Y_j - Z_j^i}{p(\theta_i)} - \frac{Y_j - Z_j^1}{p(\theta_1)} + \varepsilon \right] \qquad \text{for } i = 1, \ldots, k_\varepsilon.$$

Let $N^*(a)$ be the smallest $n$ such that, simultaneously,

$$S_n \geq a \left( 1 + \frac{\varepsilon}{p_0} \right) + \tau_\varepsilon \quad \text{and} \quad \min_{1 \leq i \leq k_\varepsilon} B_n^i \geq 0.$$

Clearly, $N^*(a) \geq N_3(a)$. Now it suffices to show that

$$(3.12) \qquad \mathbf{E}_{\lambda_0} N^*(a) \leq (1 + r_\varepsilon) \frac{a}{q(\lambda_0)}$$

for all sufficiently large $a$ for some $r_\varepsilon > 0$ which can be made arbitrarily small by choosing a sufficiently small $\varepsilon$.

To prove (3.12), assume that $\{U_i\}$ are indexed (re-index if necessary) so that the minimum (over $i$) of the left-hand side of (3.11) occurs when $i = 1$. By the proof of Lemma 2 in [7], we have

$$(3.13) \qquad \mathbf{E}_{\lambda_0} N^*(a) \leq \mathbf{E}_{\lambda_0}(v_1) + \mathbf{E}_{\lambda_0}(v_+) \mathbf{E}_{\lambda_0}(w),$$

where

$$v_1 = \inf \left\{ n : S_n > a \left( 1 + \frac{\varepsilon}{p_0} \right) + \tau_\varepsilon \right\},$$

$$v_+ = \inf \{ n : S_n > 0 \},$$

$$w = last \text{ time } \min_{1 \leq i \leq k_\varepsilon} B_n^i < 0.$$

By (3.11) and the definition of $q(\lambda)$,

$$\mathbf{E}_{\lambda_0} \left[ \frac{Y_j - Z_j^1}{p(\theta_1)} - \varepsilon \right] \geq q(\lambda_0) - \varepsilon \left( 1 + \frac{1}{p_0} \right).$$

Thus, if we choose $\varepsilon$ small enough so that $q(\lambda_0) - \varepsilon(1 + 1/p_0) > 0$, then it is well known from renewal theory that

$$\mathbf{E}_{\lambda_0}(v_1) \leq (1 + o(1)) \frac{a(1 + \varepsilon/p_0) + \tau_\varepsilon}{q(\lambda_0) - \varepsilon(1 + 1/p_0)} \quad \text{and} \quad \mathbf{E}_{\lambda_0}(v_+) = D(\varepsilon) < \infty.$$

Moreover, $\mathbf{E}_{\lambda_0}(w) = h(\varepsilon) < \infty$ because the summands in $B_n^i$ have positive mean and finite variance under $\mathbf{P}_{\lambda_0}$; see, for example, Theorem D in [7]. Relation (3.12) follows at once from (3.13). Therefore, the lemma holds. $\square$

PROOF OF THEOREM 3.1 AND COROLLARY 3.1. First we establish an upper bound of $\log \mathbf{E}_\theta T^*(a)$. By Lemma 3.2 and Lorden's lower bound (1.2),

$$\log \mathbf{E}_\theta T^*(a) \leq \inf_{\lambda \in \Lambda} ((1 + o(1)) I(\lambda, \theta) \overline{\mathbf{E}}_\lambda T^*(a)) \leq \inf_{\lambda \in \Lambda} \left( (1 + o(1)) I(\lambda, \theta) \frac{a}{q(\lambda)} \right).$$

The compactness of $\Lambda$ leads to

$$\log \mathbf{E}_\theta T^*(a) \leq (1 + o(1)) \left( \inf_{\lambda \in \Lambda} \frac{I(\lambda, \theta)}{q(\lambda)} \right) a.$$

If $p(\theta)$ is an optimizer, then $(p(\theta), q(\lambda))$ is an optimizer pair by Proposition 3.1. Thus

$$\log \mathbf{E}_\theta T^*(a) \leq (1 + o(1)) p(\theta) a.$$

Combining this with Lemma 3.1 yields

$$\log \mathbf{E}_\theta T^*(a) \sim p(\theta) a.$$

Similarly,

$$\overline{\mathbf{E}}_\lambda T^*(a) \sim a/q(\lambda),$$

and the same results are true if $T^*(a)$ is replaced by $T_1^*(a)$.

To prove Theorem 3.1, note that the asymptotic efficiency of $T^*(a)$ and $T_1^*(a)$ at $(\theta, \lambda)$ is

$$e(\theta, \lambda) = \frac{p(\theta) q(\lambda)}{I(\lambda, \theta)},$$

and so they are asymptotically optimal to first order by virtue of the compactness of $\theta$ and $\Lambda$ and the definition of an optimizer pair.

Applying Lorden's lower bound, we can prove Corollary 3.1 in the same way as the upper bound for $\log \mathbf{E}_\theta T^*(a)$. $\square$

3.2. *Optimizer pairs.* The following are some examples of an optimizer pair $(p, q)$ and the corresponding asymptotically optimal procedures.

EXAMPLE 3.1. If there exists $I_0$ such that for all $\theta \in \Theta$, $\inf_{\lambda \in \Lambda} I(\lambda, \theta) = I_0$, then $q_0(\lambda) \equiv I_0$ yields

$$p(\theta) = 1 \quad \text{and} \quad q(\lambda) = \inf_{\theta \in \Theta} I(\lambda, \theta).$$

This is even true for composite $\Theta$ and $\Lambda$. In particular, if $\Theta$ is simple, say $\{\theta_0\}$, then our consideration reduces to the standard formulation where the pre-change distribution is completely specified. Moreover, Pollak [18] proved that $T(a)$, defined in (3.4), has a *second-order* optimality property in the context of open-ended hypothesis testing if $f_\theta$ and $g_\lambda$ belong to exponential families.

EXAMPLE 3.2. If there exists $I_0$ such that for all $\lambda \in \Lambda$, $\inf_{\theta \in \Theta} I(\lambda, \theta) = I_0$, then $q_0(\lambda) \equiv 1$ yields

$$p(\theta) = \inf_{\lambda \in \Lambda} I(\lambda, \theta) \quad \text{and} \quad q(\lambda) = 1,$$

even for composite $\Theta$ and $\Lambda$. In particular, if $\Lambda$ is simple, say $\{\lambda\}$, then the considerations of Section 3 reduce to those of the problem in Section 2.

EXAMPLE 3.3. Suppose $f_\theta$ and $g_\lambda$ are exponentially distributed with unknown means $1/\theta$ and $1/\lambda$, respectively. Assume $\Theta = \{\theta : \theta \in [\theta_0, \theta_1]\}$ and $\Lambda = \{\lambda : \lambda \in [\lambda_0, \lambda_1]\}$, where $\theta_0 < \theta_1 < \lambda_0 < \lambda_1$. Then optimizer pairs $(p(\theta), q(\lambda))$ are not unique. For example, the following two pairs are nonequivalent:

$$\begin{cases} p_1(\theta) = I(\lambda_0, \theta), \\ q_1(\lambda) = I(\lambda, \theta_0)/I(\lambda_0, \theta_0), \end{cases} \quad \text{and} \quad \begin{cases} p_2(\theta) = I(\lambda_1, \theta)I(\lambda_0, \theta_1)/I(\lambda_1, \theta_1), \\ q_2(\lambda) = I(\lambda, \theta_1)/I(\lambda_0, \theta_1). \end{cases}$$

Suppose $t_1^*(a)$ and $t_2^*(a)$ are the procedures defined by (3.5) for the pairs $(p_1(\theta), q_1(\lambda))$ and $(p_2(\theta), q_2(\lambda))$, respectively. Even though both $t_1^*(a)$ and $t_2^*(a)$ are asymptotically optimal to first order, $t_1^*(a)$ performs better uniformly over $\Theta$ (in the sense of larger long ARL), while $t_2^*(a)$ performs better uniformly over $\Lambda$ (in the sense of smaller short ARL).

3.3. *Numerical simulations.* In this section we report some simulation studies comparing the performance of our procedures with a commonly used procedure in the literature.

The simulations consider the problem of detecting a change in distribution from $f_\theta$ to $g_\lambda$, where $f_\theta$ and $g_\lambda$ are exponentially distributed with unknown means $1/\theta$ and $1/\lambda$, respectively, and $\theta \in \Theta = [0.8, 1]$ and $\lambda \in \Lambda = [2, 3]$.

Note that $q_0(\lambda) \equiv 1$ leads to an optimizer $p(\theta) = I(2, \theta)$ where $I(\lambda, \theta) = \theta/\lambda - 1 - \log(\theta/\lambda)$, and so our procedure based on (3.6) is defined by

$$\hat{T}^*(a) = \inf\left\{ n \geq 1 : \max_{1 \leq k \leq n} \inf_{0.8 \leq \theta \leq 1} \sup_{2 \leq \lambda \leq 3} \frac{\lambda - \theta}{p(\theta)} \sum_{i=k}^{n} \left( \frac{\log \lambda - \log \theta}{\lambda - \theta} - X_i \right) \geq a \right\}.$$

A commonly used procedure in the change-point literature is the *generalized likelihood ratio procedure* which specifies the nominal value $\theta_0$ (of the parameter of the pre-change distribution); see [14] and [29]. The procedure is defined by the stopping time

$$\tau(\theta_0, a) = \inf\left\{ n \geq 1 : \max_{1 \leq k \leq n} \sup_{\lambda \in \Lambda} \sum_{i=k}^{n} \log \frac{g_\lambda(X_i)}{f_{\theta_0}(X_i)} \geq a \right\}$$

$$= \inf\left\{ n \geq 1 : \max_{1 \leq k \leq n} \sup_{2 \leq \lambda \leq 3} \sum_{i=k}^{n} \left( \log \frac{\lambda}{\theta_0} - (\lambda - \theta_0)X_i \right) \geq a \right\}.$$

Note that $\tau(\theta_0, a)$ can be thought of as our procedure $\hat{T}^*(a)$ whose $\Theta$ contains the single point $\theta_0$. The choice of $\theta_0$ can be made directly by considering the pre-change distribution which is *closest* to the post-change distributions because it is always more difficult to detect a smaller change. For our example, $\theta_0 = 1$.

An effective method to implement $\tau(\theta_0, a)$ numerically can be found in [14]. Similarly, we can implement $\hat{T}^*(a)$ as follows. Compute $V_n$ recursively by $V_n = \max(V_{n-1} + \log(2/0.8) - (2 - 0.8)X_n, 0)$. Whenever $V_n = 0$, one can begin a new cycle, discarding all previous observations and starting afresh on the incoming observations, because for all $0.8 \leq \theta \leq 1$, $2 \leq \lambda \leq 3$ and $1 \leq k \leq n$, $\sum_{i=k}^{n}((\log \lambda - \log \theta)/(\lambda - \theta) - X_i) \leq 0$ since $(\log \lambda - \log \theta)/(\lambda - \theta)$ is maximized at $(\theta, \lambda) = (0.8, 2)$. Now each time a new cycle begins compute at each stage $n = 1, 2, \ldots$

$$Q_k^{(n)} = X_n + \cdots + X_{n-k+1}, \qquad k = 1, \ldots, n.$$

Then the procedure $\hat{T}^*(a) =$ first $n$ such that $Q_k^{(n)} < c_k$ for some $k$, where

$$c_k = \inf_{0.8 \leq \theta \leq 1} \sup_{2 \leq \lambda \leq 3} \left[ k \frac{\log \lambda - \log \theta}{\lambda - \theta} - \frac{p(\theta)a}{\lambda - \theta} \right].$$

To further speed up the implementation, compute $W_n$ recursively by $W_n = \max(W_{n-1} + \log 2 - X_n, 0)$. Stop whenever $W_n \geq p(0.8)a/1.2$. Continue taking new observations (i.e., do not stop) whenever $W_n \leq p(1)a/2$. If $p(1)a/2 \leq W_n \leq p(0.8)a/1.2$, then we will also stop at time $n$ if $Q_k^{(n)} < c_k$ for some $k$. The reasons behind this implementation are given below.

First, if at time $n_0$ we have $W_{n_0} \geq p(0.8)a/1.2 > 0$, then there exists some $k_0$ such that $\sum_{i=k_0}^{n_0}(\log 2 - X_i) \geq p(0.8)a/1.2$. Thus for all $\theta \in [0.8, 1]$ and $\lambda_0 = 2$,

$$\frac{\lambda_0 - \theta}{p(\theta)} \sum_{i=k_0}^{n_0} \left( \frac{\log \lambda_0 - \log \theta}{\lambda_0 - \theta} - X_i \right) \geq \frac{\lambda_0 - \theta}{p(\theta)} \sum_{i=k_0}^{n_0} (\log 2 - X_i)$$

$$\geq \frac{\lambda_0 - \theta}{p(\theta)} \cdot \frac{p(0.8)a}{1.2} \geq a.$$

*Comparison of two procedures in change-point problems with composite pre-change and composite post-change hypotheses*

|  | $a$ | $\hat{T}^*(a)$ 22.50 | $\tau(1,a)$ 5.02 |
|---|---|---|---|
| $\mathbf{E}_\theta N$ | $\theta = 1$ | $601 \pm 18$ | $606 \pm 19$ |
|  | $\theta = 0.9$ | $1448 \pm 43$ | $1207 \pm 36$ |
|  | $\theta = 0.8$ | $3772 \pm 116$ | $2749 \pm 90$ |
| $\mathbf{E}_\lambda N$ | $\lambda = 2$ | $21.41 \pm 0.10$ | $21.92 \pm 0.11$ |
|  | $\lambda = 2.2$ | $18.09 \pm 0.07$ | $18.18 \pm 0.09$ |
|  | $\lambda = 2.5$ | $15.08 \pm 0.05$ | $14.76 \pm 0.06$ |
|  | $\lambda = 2.7$ | $13.75 \pm 0.04$ | $13.22 \pm 0.05$ |
|  | $\lambda = 3$ | $12.29 \pm 0.04$ | $11.62 \pm 0.04$ |

Hence, $\hat{T}^*(a)$ will stop at time $n_0$. Second, $\hat{T}^*(a)$ will never stop at time $n$ when $W_n \leq p(1)a/2$ because for $\theta_1 = 1$, all $2 \leq \lambda \leq 3$, and all $k$,

$$\frac{\lambda - \theta_1}{p(\theta_1)} \sum_{i=k}^{n} \left( \frac{\log \lambda - \log \theta_1}{\lambda - \theta_1} - X_i \right) \leq \frac{\lambda - \theta_1}{p(\theta_1)} \sum_{i=k}^{n} (\log 2 - X_i) \leq \frac{\lambda - \theta_1}{p(\theta_1)} W_n \leq a.$$

Table 2 provides a comparison of the performances for our procedure $\hat{T}^*(a)$ with those of $\tau(\theta_0, a)$. The threshold $a$ for each of these two procedures is determined from the criterion $\mathbf{E}_{\theta=1} N(a) \approx 600$. The results in Table 2 are based on 1000 simulations for $\mathbf{E}_\theta N$ and 10,000 simulations for $\mathbf{E}_\lambda N$. Note that for these two procedures, the detection delay $\overline{\mathbf{E}}_\lambda N = \mathbf{E}_\lambda N$. Table 2 shows that at a small additional cost of detection delay, $\hat{T}^*(a)$ can significantly improve the mean times between false alarms compared to $\tau(1,a)$. This is consistent with the asymptotic theory in this section.

**4. Normal distributions.** Our general theory in Section 3 assumes that $\Theta$ and $\Lambda$ are compact. If they are not compact, then our proposed procedures *may* or *may not* be asymptotically optimal. However, we can still sometimes apply our ideas in these situations, as shown in the following example.

Suppose we want to detect a change from negative to positive in the mean of independent normally distributed random variables with variance 1. In the context of open-ended hypothesis testing, we want to test

$$H_0 : \theta \in \Theta = (-\infty, 0) \quad \text{against} \quad H_1 : \lambda \in \Lambda = (0, \infty).$$

Let us examine the procedures $\hat{T}(a)$ defined in (3.6) for different choices of optimizer pairs.

First, let us assume $q_0(\lambda) = \lambda^{1/\beta}$ with $\beta \geq 1/2$; then we have an optimizer pair

$$p(\theta) = k_\beta |\theta|^{2 - (1/\beta)} \quad \text{and} \quad q(\lambda) = \lambda^{1/\beta} \qquad \text{with } k_\beta = 2\beta^2 (2\beta - 1)^{(1/\beta) - 2}$$

(assume $0^0 = 1$), and thus the procedure defined in (3.6) becomes $\hat{t}_\beta(a) =$ first time $n$ such that

$$\inf_{\theta < 0} \sup_{\lambda > 0} \left[ \frac{1}{p(\theta)} \left( (\lambda - \theta) S_n - \frac{\lambda^2 - \theta^2}{2} n \right) \right] \geq a \qquad \text{where } S_n = \sum_{i=1}^{n} X_i.$$

Letting $\theta \to 0$ gives us that $S_n > 0$ if $\hat{t}_\beta(a) = n$, and rewriting the stopping rule as

$$\inf_{\theta < 0} \sup_{\lambda > 0} \left[ -\left( \lambda - \frac{S_n}{n} \right)^2 + \left( \frac{S_n}{n} - \theta \right)^2 - \frac{2}{n} p(\theta) a \right] \geq 0.$$

The supremum is attained at $\lambda = S_n/n$, and so $\hat{t}_\beta(a) =$ first time $n$ such that for all $\theta < 0$,

$$\frac{S_n}{n} \geq \theta + \sqrt{\frac{2}{n} p(\theta) a}.$$

A routine calculation leads to

$$\hat{t}_\beta(a) = \inf\{n \geq 1 : S_n \geq a^\beta n^{1-\beta}\}.$$

This suggests using a stopping time of the form

(4.1) $$\hat{t}_\beta^*(a) = \inf\left\{ n \geq 1 : \max_{0 \leq k \leq n} [(S_n - S_k)(n-k)^{\beta-1}] \geq a^\beta \right\}$$

to detect a change in mean from negative to positive. Observe that for $\beta = 1$, $\hat{t}_\beta(a)$ is just the one-sided SPRT and $\hat{t}_\beta^*(a)$ is just a special form of Page's CUSUM procedures. For $\beta = 1/2$, $\hat{t}_\beta(a)$ and $\hat{t}_\beta^*(a)$ have also been studied extensively in the literature, since they are based on the generalized likelihood ratio. Different motivation to obtain these two procedures can be found for $\hat{t}_\beta(a)$ in Chapter IV of [28], which is from the viewpoint of the repeated significant test, and for $\hat{t}_\beta^*(a)$ in [29], which is from the viewpoint of the generalized likelihood ratio. For $\hat{t}_\beta(a)$ with $0 < \beta \leq 1$, see [3] and equation (9.2) on page 188 in [28].

Next, $q_0(\lambda) = 1$ leads to

$$p(\theta) = \frac{\theta^2}{2} \quad \text{and} \quad q(\lambda) = 1$$

and

$$\hat{t}_0(a) = \inf\{n \geq a : S_n \geq 0\}.$$

Hence we use the following stopping time to detect a change in mean from negative to positive:

(4.2) $$\hat{t}_0^*(a) = \inf\left\{ n \geq a : \max_{0 \leq k \leq n-a} (S_n - S_k) \geq 0 \right\},$$

where the maximum is taken over $0 \leq k \leq n - a$. It is interesting to see that $\hat{t}_0(a)$ and $\hat{t}_0^*(a)$ can be thought of as the limits of $\hat{t}_\beta(a)$ and $\hat{t}_\beta^*(a)$, respectively, as $\beta \to \infty$.

Though one cannot use our theorems directly to analyze the properties of $\hat{t}_0^*(a)$ and $\hat{t}_\beta^*(a)$, they are indeed asymptotically optimal to first order. For $\beta \geq 1/2$, first note that

$$p(\theta)q(\lambda) = I(\lambda, \theta) \qquad \text{if } \theta = -(2\beta - 1)\lambda.$$

By nonlinear renewal theory ([28], Chapters 9 and 10),

$$\overline{\mathbf{E}}_\lambda \hat{t}_\beta^*(a) \sim a/q(\lambda).$$

Equation (13) on page 636 in [15] shows that for any $\theta \leq 0$,

$$\mathbf{P}_\theta(\hat{t}_\beta(a) < \infty) \leq \exp(-(1 + o(1))p(\theta)a),$$

and so Lemma 2.1 implies $\log \mathbf{E}_\theta \hat{t}_\beta^*(a) \sim p(\theta)a$ as $a \to \infty$. Thus $\hat{t}_\beta^*(a)$ is asymptotically efficient at $(\theta, \lambda)$ with $\theta = -(2\beta - 1)\lambda$, and hence $\hat{t}_\beta^*(a)$ is asymptotically optimal to first order. Similarly, the asymptotic optimality property of $\hat{t}_0^*(a)$ can be proved directly since the structure of $\hat{t}_0(a)$ is very simple.

REMARK. The above arguments establish the following optimality properties of $\hat{t}_\beta(a)$ and $\hat{t}_\beta^*(a)$. Suppose we want to test

$$H_{0,\delta} : \theta = -(2\beta - 1)\delta \quad \text{against} \quad H_{1,\delta} : \lambda = \delta,$$

where $\beta \geq 1/2$ is given but $\delta > 0$ is unknown. Then $\hat{t}_\beta(a)$ is an asymptotically optimal solution for all $\delta > 0$, while $\hat{t}_\beta^*(a)$ is asymptotically optimal in the problems of detecting a change from $H_{0,\delta}$ to $H_{1,\delta}$ for all $\delta > 0$. As far as we know, no optimality properties of $\hat{t}_\beta(a)$ and $\hat{t}_\beta^*(a)$ have been studied except for the special case of $\beta = 1/2$ or 1. Even for the case $\beta = 1/2$ which was studied in [29], our method is simpler and more instructive.

**5. Proof of Theorem 2.1.** The basic idea in proving Theorem 2.1 is to relate the stopping time $M(a)$ in (2.4) to new stopping times defined by

$$(5.1) \qquad M_\theta(a) = \inf\left\{ n \geq 1 : \sum_{i=1}^{n} \log \frac{f_\lambda(X_i)}{f_\theta(X_i)} - I(\lambda, \theta)a > 0 \right\}.$$

The proof of Theorem 2.1 is based on the following lemmas.

LEMMA 5.1. *For all $\theta \in [\theta_0, \theta_1]$,*

$$\mathbf{P}_\theta(M(a) < \infty) \leq \mathbf{P}_\theta(M_\theta(a) < \infty) \leq \exp(-I(\lambda, \theta)a),$$

*and hence* (2.7) *holds.*

PROOF.   The first inequality follows at once from the fact that $M(a) \geq M_\theta(a)$ for all $\theta \in [\theta_0, \theta_1]$, and the second inequality is a direct application of Wald's likelihood ratio identity.   □

We now derive approximations for $\mathbf{E}_\lambda M(a)$. Similarly to (2.6), $M_\theta(a)$ in (5.1) can be written as

$$M_\theta(a) = \inf\{n \geq 1 : S_n \geq b'(\lambda)a + (n-a)\phi(\theta)\}.$$

As we said earlier, the supremum in (2.6) is attained at $\theta = \theta_0$ if $n \leq a$, and at $\theta = \theta_1$ if $n > a$, so that

(5.2)        $\{M(a) = m\} = \{M(a) = M_{\theta_0}(a) = m\}$        for all $m \leq a$.

For simplicity, we omit $a$ and $\theta$, writing $M = M(a)$ and $M_k = M_{\theta_k}(a)$ for $k = 0, 1$.

LEMMA 5.2.   As $a \to \infty$,

$$\mathbf{E}_\lambda M(a) = a + \frac{\phi(\theta_1) - \phi(\theta_0)}{b'(\lambda) - \phi(\theta_1)} \mathbf{E}_\lambda(a - M_0; M_0 \leq a) + O(1).$$

PROOF.   Observe that

$$\mathbf{E}_\lambda M = a - \mathbf{E}_\lambda(a - M; M \leq a) + \mathbf{E}_\lambda(M - a; M > a),$$

and by (5.2), $\mathbf{E}_\lambda(M - a; M \leq a) = \mathbf{E}_\lambda(M_0 - a; M_0 \leq a)$. Thus it suffices to show that

(5.3)   $\mathbf{E}_\lambda(M - a; M > a) = \dfrac{b'(\lambda) - \phi(\theta_0)}{b'(\lambda) - \phi(\theta_1)} \mathbf{E}_\lambda(a - M_0; M_0 \leq a) + O(1).$

To prove this, define a stopping time

$$N_k(u) = \inf\left\{ n \geq 1 : \sum_{i=1}^{n}(X_i - \phi(\theta_k)) \geq u \right\},$$

for $k = 0, 1$ and any $u > 0$. Assume $a$ is an integer. (For general $a$, using $[a]$, the largest integer $\leq a$, permits one to carry through the following argument with minor modifications.) By (5.2) we have

$$\mathbf{E}_\lambda(M - a | M > a) = \int_{-\infty}^{0} \mathbf{E}_\lambda(M - a | S_a - b'(\lambda)a = x, M_0 > a)$$
$$\times \mathbf{P}_\lambda(S_a - b'(\lambda)a \in dx | M_0 > a).$$

Conditional on the event $\{S_a - b'(\lambda)a = x, M_0 > a\}$,

$$M - a = \inf\{m : X_{a+1} + \cdots + X_{a+m} + S_a \geq b'(\lambda)a + m\phi(\theta_1)\}$$
$$= \inf\left\{ m : \sum_{i=1}^{m}(X_{a+i} - \phi(\theta_1)) \geq b'(\lambda)a - S_a = -x \right\},$$

which is equivalent to $N_1(-x)$ since $X_1, X_2, \ldots$ are independent and identically distributed. Thus

(5.4)
$$\mathbf{E}_\lambda(M - a | M > a)$$
$$= \int_{-\infty}^0 \mathbf{E}_\lambda N_1(-x) \mathbf{P}_\lambda(S_a - b'(\lambda)a \in dx | M_0 > a).$$

Similarly,

(5.5)
$$\mathbf{E}_\lambda(M_0 - a | M_0 > a)$$
$$= \int_{-\infty}^0 \mathbf{E}_\lambda N_0(-x) \mathbf{P}_\lambda(S_a - b'(\lambda)a \in dx | M_0 > a).$$

Now for $k = 0, 1$ and any $u > 0$, define

$$R_k(u) = \sum_{i=1}^{N_k(u)} (X_i - \phi(\theta_k)) - u.$$

Then, by Theorem 1 in [13],

$$\sup_{u \geq 0} \mathbf{E}_\lambda R_k(u) \leq \mathbf{E}_\lambda (X_1 - \phi(\theta_k))^2 / (b'(\lambda) - \phi(\theta_k)) < \infty.$$

By Wald's equation, $(b'(\lambda) - \phi(\theta_k)) \mathbf{E}_\lambda N_k(u) = u + \mathbf{E}_\lambda R_k(u)$, so that

$$\sup_{u \geq 0} \mathbf{E}_\lambda \left( N_k(u) - \frac{u}{b'(\lambda) - \phi(\theta_k)} \right) < \infty$$

for $k = 0, 1$. Hence, we have

$$\sup_{u \geq 0} \left| \mathbf{E}_\lambda N_1(u) - \frac{b'(\lambda) - \phi(\theta_0)}{b'(\lambda) - \phi(\theta_1)} \mathbf{E}_\lambda N_0(u) \right| < \infty.$$

Plugging into (5.4), and comparing with (5.5), we have

$$\mathbf{E}_\lambda(M - a | M > a) = \frac{b'(\lambda) - \phi(\theta_0)}{b'(\lambda) - \phi(\theta_1)} \mathbf{E}_\lambda(M_0 - a | M_0 > a) + O(1).$$

Relation (5.3) follows at once from the fact that $\{M > a\} = \{M_0 > a\}$ and the fact that $\mathbf{E}_\lambda(M_0 - a) = O(1)$. Hence, the lemma holds. $\square$

LEMMA 5.3. *Suppose $Y_1, Y_2, \ldots$ are independent and identically distributed with mean $\mu > 0$ and finite variance $\sigma^2$. Define*

$$N_a = \inf \left\{ n \geq 1 : \sum_{i=1}^n Y_i \geq a \right\}.$$

*Then as $a \to \infty$,*

$$\mathbf{E} \left( \frac{a}{\mu} - N_a ; N_a \leq \frac{a}{\mu} \right) = \sqrt{a} \left( \frac{\sigma}{\sqrt{2\pi \mu^3}} + o(1) \right).$$

PROOF.   The lemma follows at once from the well-known facts that as $a \to \infty$,

$$\mathbf{E}(N_a) = \frac{a}{\mu} + O(1) \quad \text{and} \quad \text{Var}(N_a) = (1 + o(1))\frac{\sigma^2 a}{\mu^3},$$

and that

$$\frac{N_a - a/\mu}{\sqrt{a\sigma^2/\mu^3}}$$

is asymptotically standard normal. See page 372 in [4], equation (5) in [27] and Theorem 8.34 in [28].   □

PROOF OF THEOREM 2.1.   Relation (2.7) is proved in Lemma 5.1. By (5.1) $M_0 = M_{\theta_0}(a)$ can be written as

$$M_0 = \inf\left\{n \geq 1 : \sum_{i=1}^{n} \frac{1}{I(\lambda, \theta_0)} \log \frac{f_\lambda(X_i)}{f_{\theta_0}(X_i)} > a\right\}.$$

By Lemma 5.3 it is easy to show that

$$\mathbf{E}_\lambda(a - M_0; M_0 \leq a) = \sqrt{a}\left(\frac{\sigma_0}{\sqrt{2\pi}} + o(1)\right),$$

where $\sigma_0 = \sqrt{b''(\lambda)}/(b'(\lambda) - \phi(\theta_0))$. Thus relation (2.8) holds by Lemma 5.2 and the definition of $\phi(\theta)$ in (2.5).   □

**Acknowledgments.**   This work is part of my Ph.D. dissertation at the California Institute of Technology. I would like to thank my thesis advisor, Professor Gary Lorden, for his constant support and encouragement and Professor Moshe Pollak for sharing his insightful ideas. Thanks also to Professor Sarah Holte, Professor Jon A. Wellner, the Associate Editor and the referees for their helpful remarks.

FRED HUTCHINSON CANCER RESEARCH CENTER
1100 FAIRVIEW AVENUE NORTH, M2-B500
SEATTLE, WASHINGTON 98109
USA
E-MAIL: ymei@fhcrc.org